\documentclass[a4paper]{article}

\usepackage[english]{babel}
\usepackage[utf8x]{inputenc}
\usepackage[T1]{fontenc}

\usepackage[a4paper,top=2cm,bottom=2.5cm,left=2.3cm,right=2.3cm,marginparwidth=1.75cm,margin=1.0in]{geometry}

\usepackage[colorinlistoftodos]{todonotes}
\usepackage[colorlinks=true, allcolors=blue]{hyperref}
\usepackage{xy}
\usepackage[affil-it]{authblk}
\usepackage{amssymb}
\usepackage{cancel}
\usepackage{epstopdf}
\usepackage{verbatim}
\usepackage{amsmath,amscd}
\usepackage{graphicx}
\usepackage{subcaption}
\usepackage{mathtools}
\usepackage{amsthm}

\usepackage{mathrsfs}
\usepackage{tikz}
\usetikzlibrary{cd}
\usepackage{pgf}
\usepackage{mathtools}

\usepackage[noend]{algpseudocode}
\usepackage{algorithm}

\usepackage{mathrsfs}
\usetikzlibrary{arrows}
\usepackage[toc,page]{appendix}
\usepackage{caption}
\usepackage[colorinlistoftodos]{todonotes}

\graphicspath{ {./Images/} }

\makeatletter
\newcommand{\subjclass}[2][2010]{
  \let\@oldtitle\@title
  \gdef\@title{\@oldtitle\footnotetext{#1 {Mathematics subject classification:} #2}}
}

\newcommand\sponsor[1]{%
  \begingroup
  \renewcommand\thefootnote{}\footnote{#1}%
  \addtocounter{footnote}{-1}%
  \endgroup
}

\newtheorem{theorem}{Theorem}[section]

\makeatletter
\newtheorem*{rep@theorem}{\rep@title}
\newcommand{\newreptheorem}[2]{
\newenvironment{rep#1}[1]{
 \def\rep@title{#2 \ref{##1}}
 \begin{rep@theorem}}
 {\end{rep@theorem}}}
\makeatother

 \newreptheorem{theorem}{Theorem}

\theoremstyle{definition}
\newtheorem{mydef}[theorem]{Definition}

\newtheorem{rem}[theorem]{Remark}
\newtheorem{exa}[theorem]{Example}
\newtheorem{cor}[theorem]{Corollary}

\title{{\bf{Tropical Igusa Invariants}}} 
\author{\large{Paul Alexander Helminck}}
\affil{Durham University\\ 
\vspace{0.3cm}
Department of Mathematics}


\date{}
\begin{document}
\maketitle
\definecolor{qqqqff}{rgb}{0,0,1}
\begin{abstract}
Let $X$ be a smooth geometrically connected projective curve of genus two over a complete non-archimedean field $K$. For discretely valued $K$, 
 the first main theorem in \cite{liu} gives a set of criteria on the Igusa invariants of the curve that determine the minimal skeleton of $X$ together with its edge lengths and vertex weights. In this paper we use the theory of Berkovich spaces to give a new proof of this theorem that works for arbitrary complete non-archimedean fields. 
 We furthermore interpret the final result in terms of tropical moduli spaces and tropical Igusa invariants. This reformulation shows that the abstract tropicalization map ${M}_{2}\to\mathrm{trop}(M_{2})$ factors through the tropicalization of a concrete embedding of ${M}_{2}$ into a weighted projective space. 

\end{abstract}

\sponsor{The author was supported by the UKRI Fellowship {\it{"Computational tropical geometry and its applications"}} with reference number MR/S034463/1.}

\section{Introduction}\label{IntroductionSection}

\begin{figure}[h!]
\centering
\includegraphics[height=6.5cm]{{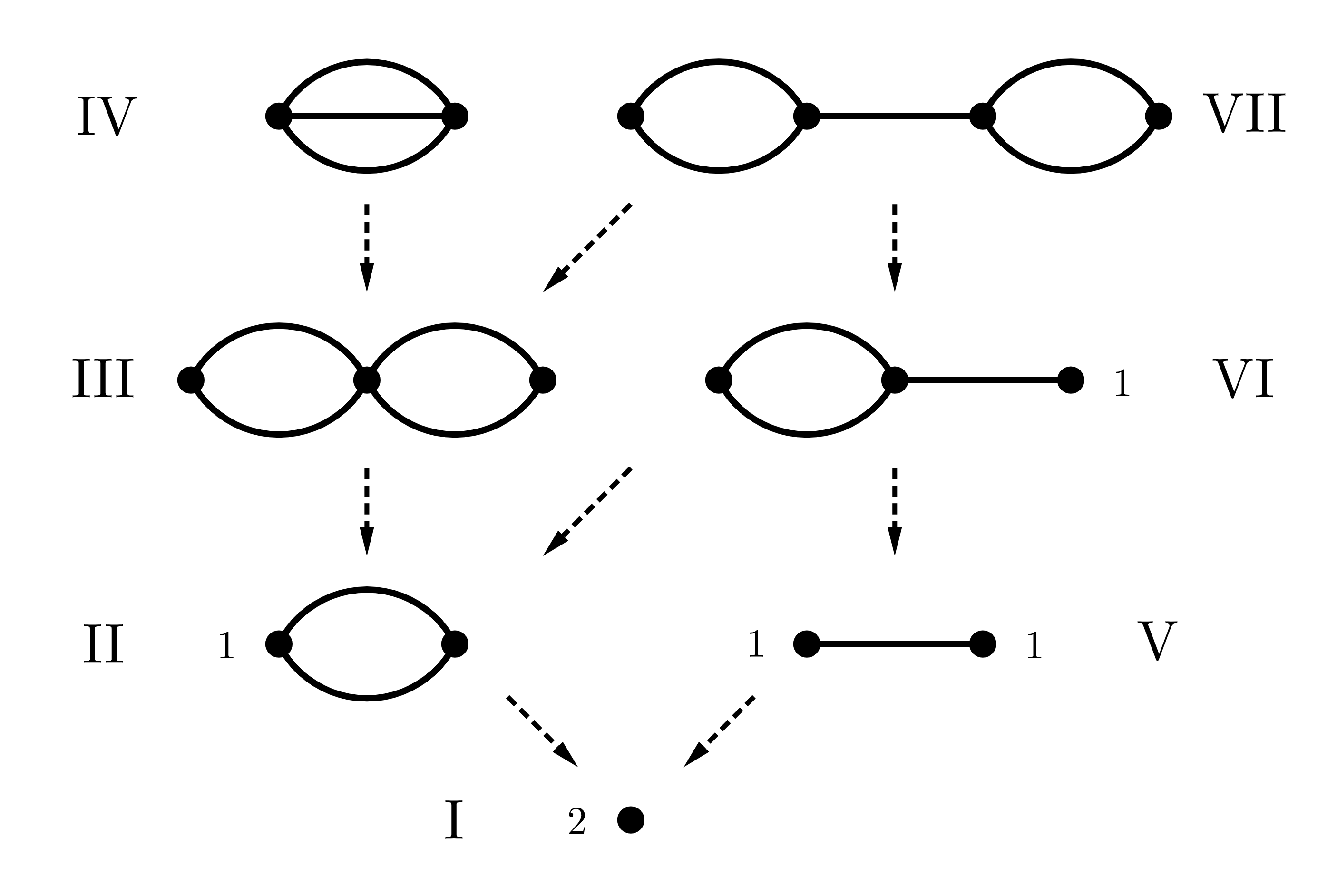}}
\caption{\label{Genus2Mod}
The moduli space of tropical curves of genus $2$. The Roman numerals denote the reduction type and the Arabic numerals denote the genera of the vertices.  
} 
\end{figure}

Let $X/K$ be a smooth proper curve of genus two over a complete algebraically closed non-archimedean field $K$. An important invariant of the Berkovich analytification of $X$ is its {{minimal skeleton}} $\Sigma(X)$, which can be described as the set of points in $X^{\mathrm{an}}$ that do not admit an affinoid neighborhood isomorphic to a closed disk. Alternatively, it is the dual graph of the special fiber of a stable model of $X$ together with some additional data on its edges and vertices.
For curves of genus two, there are seven possible graph-theoretical types for $\Sigma(X)$, see Figure \ref{Genus2Mod}. 

Our goal in this paper is to characterize the minimal skeleton in terms of the Igusa invariants of $X$. To that end, we introduce a set of tropical Igusa invariants, which are the valuations of certain Igusa invariants. 
We then prove the following.

\begin{theorem}\label{MainThm1}
Suppose that the characteristic of the residue field of $K$ is not two. The tropical Igusa invariants then completely determine the minimal skeleton of a curve of genus two. The conditions that distinguish between the different reduction types are given by half-spaces. 
\end{theorem}

This result is a generalization of \cite[Th\'{e}or\`{e}me 1]{liu}, which proves the same result but over a complete discretely valued field.
Our proof also simplifies certain aspects of the proof in \cite{liu} using the framework of Berkovich spaces. For instance, we do not have to work with local affine $R$-models, since we can use the results in
\cite{InvariantsSuper} to deduce the skeleton from the tree of the branch locus of the hyperelliptic map $X\to\mathbb{P}^{1}$. The main idea after this is as follows. For every reduction type, we write down an analytic universal family. By applying suitable M\"{o}bius transformations on the tree of the branch locus, we then find that every curve with this reduction type occurs in this family. We then compute the tropical Igusa invariants for each universal family and after some further calculations using Gr\"{o}bner bases we find the statement of the theorem. 

We also connect this material to the abstract theory of tropical curves of genus two found in \cite{ACP2015}. To do this, we embed the coarse moduli space $M_{2}$ into a weighted projective space using a set of Igusa invariants. The tropicalization of this embedding forms the natural target space for the tropical Igusa invariants. We then use Theorem \ref{MainThm1} to define a surjective map from this tropicalization to the abstract tropicalization $M_{2}^{\mathrm{trop}}$. The final result can be found in Corollary \ref{CorModuliSpace}.    
 
 Before we start the main material of the paper, we give a short list of related papers. There is a similar result for elliptic curves using the $j$-invariant, which says that the minimal skeleton $\Sigma(X)$ has a cycle (of length $-v(j)$) if and only if $v(j)<0$. Various incarnations of this result can be found in the literature. For instance, on the tropical side this can be found in \cite{katz_markwig_markwig_2009}, \cite{BPR11} and \cite{FaithTropEll}. For curves of genus two, there is Liu's paper \cite{liu}, which is a precursor to the current paper. Tropically, there is 
 \cite{CuetoMarkwig2016}, which gives explicit faithful tropicalizations in some cases and a part of Theorem \ref{MainThm1} (using different half-spaces).
 
\subsection{Leitfaden}
We start in Section \ref{SectionSkeleta} by recalling the concept of a skeleton of a marked curve and the action of $\mathrm{PGL}_{2}$. We then discuss the Igusa invariants, see Section \ref{SectionIgusa}. In Section \ref{SectionTropIgusa} we tropicalize these invariants and obtain the tropical Igusa invariants. 
The formulas for the various reduction types and edge lengths of skeleta of curves of genus two can also be found in this section. We end Section \ref{SectionTropIgusa} by connecting our embedded tropicalization to the abstract tropicalizations studied in \cite{ACP2015}. We then prove Theorem \ref{MainThm1}. This theorem is split up into two parts: \ref{MainThm2} and \ref{MainThm3}. The first distinguishes between the various reduction types and the second gives the edge lengths for the skeleton of a given reduction type. Their proofs can be found in Sections \ref{ProofMainThm2} and \ref{ProofMainThm3} respectively.  

\section{The main theorem}\label{MainTheoremSection}

This paper uses many concepts and results from \cite{DecompFund1} and \cite{InvariantsSuper} on coverings of curves; the reader is referred to those papers for more details. These papers in turn heavily rely on \cite{ABBR1}, \cite{BPRa1} and \cite{BPR11}, so we recommend reading these as well. The list below gives a short summary of the notation we use in this paper. 

\begin{itemize}
\item $K$ is a complete algebraically closed non-archimedean field with valuation ring $R$, maximal ideal $\mathfrak{m}$ and residue field $R/\mathfrak{m}=k$. 
\item The valuation $K\to\mathbb{R}\cup\{\infty\}$ is denoted by $v(\cdot{})$. The absolute value associated to $v(\cdot{})$ is $e^{-v(x)}$, where $e$ is Euler's number.  
\item We fix a splitting of the valuation map and write $\varpi^{m}$ for the corresponding element of valuation $m\in\Gamma\subset{\mathbb{R}}$, where $\Gamma$ is the value group of $K$. 
\item We write $X$ or $X/K$ for a smooth irreducible proper curve over $K$. We will also simply call this a curve. 
We write $X^{\mathrm{an}}$ for its 
Berkovich analytification. 
\item Given an embedding of a variety $Y$ into a toric variety $Z$, we define its tropicalization as in \cite{Gubler2013}, \cite{Payne2009} or \cite{Rab2012}. 
\item Suppose that the Euler characteristic $\chi(X,D)=2-2g-|D|$ of a marked curve $(X,D)$ is strictly smaller than zero. The reduction type of $(X,D)$ is the weighted stable graph of its minimal skeleton 
(see Definition \ref{SkeletaDefinition} and \cite{ACP2015}).
\end{itemize}

\subsection{Skeleta of curves of genus two}\label{SectionSkeleta}

Let $X$ be a curve over $K$ and let $D\subset{X(K)}$. We assume here that the Euler characteristic of the marked curve $(X,D)$ is negative. The marked curve  
$(X,D)$ then admits a unique minimal skeleton $\Sigma(X,D)$ by \cite[Section 4.16]{BPRa1}. By adding the genera of the type-$2$ points in $X^{\mathrm{an}}$, this becomes a weighted topological graph. 
The skeleton $\Sigma(X,D)$ furthermore inherits the skeletal metric from $X^{\mathrm{an}}$, making it into a weighted metric graph. 
\begin{mydef}{\bf{[Skeleta]}}\label{SkeletaDefinition}
The weighted metric graph $\Sigma(X,D)$ associated to the marked curve $(X,D)$ is the minimal skeleton of $(X,D)$.  
\end{mydef}
We are interested in this minimal skeleton for curves of genus two with $D=\emptyset$. The total genus function
 \begin{equation}
 g(\Sigma(X))=\beta_{1}(\Sigma(X))+\sum_{v\in\Sigma(X)}g(v)
 \end{equation}
 is then equal to two. Here $\beta_{1}(\Sigma(X))$ is the first Betti number of $\Sigma(X)$ and $g(v)$ is the genus of the curve associated to $v$. This gives seven different reduction types, which can be found in Figure \ref{Genus2Mod}. The arrows in this figure will be discussed later.

\begin{figure}[h!]
\centering
\includegraphics[height=6.5cm]{{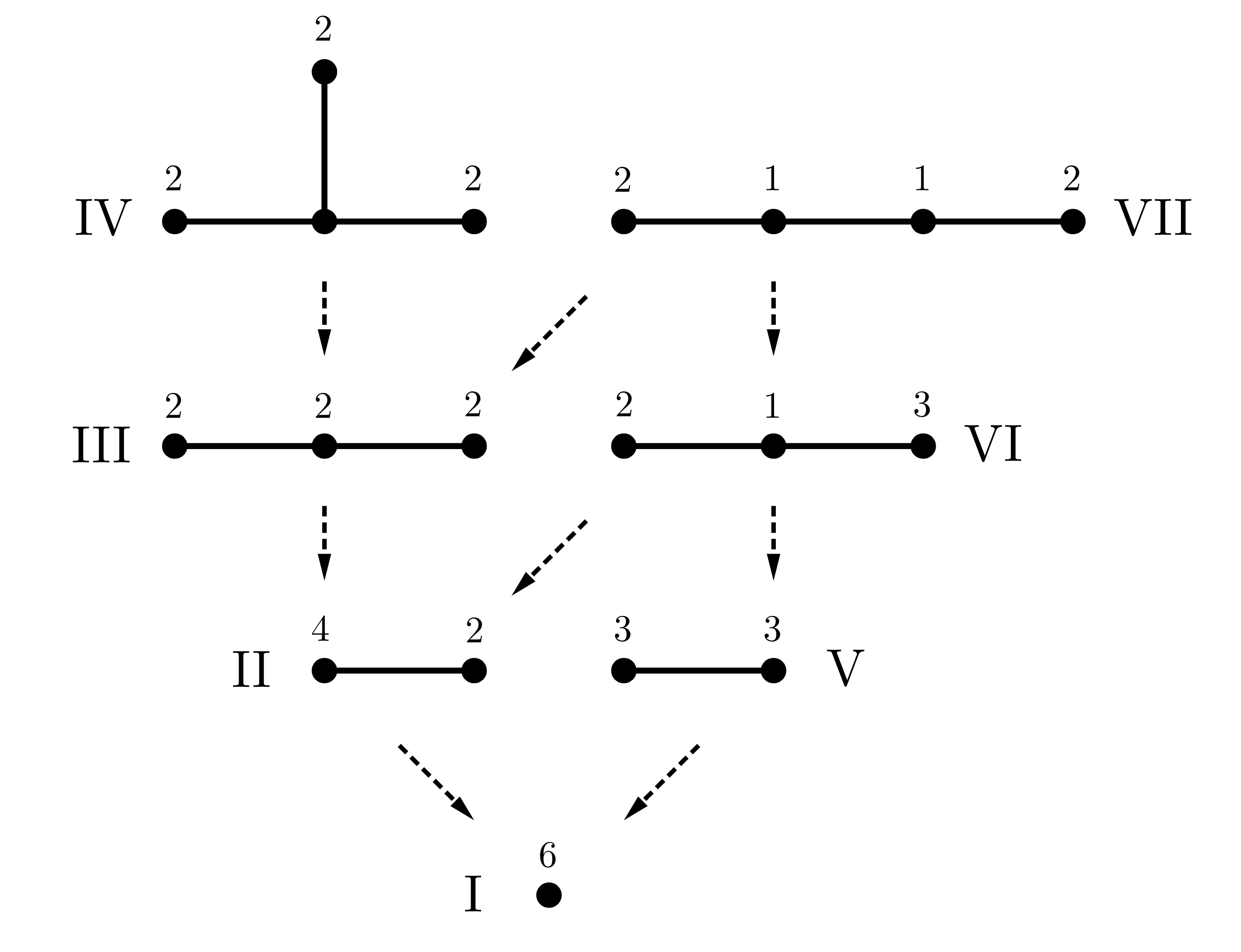}}
\caption{\label{TreesSixLeaves}
Phylogenetic trees on six leaves. 
The Roman numerals denote the type and the Arabic numerals denote the number of leaves attached to a vertex.  
} 
\end{figure}

The material in \cite{InvariantsSuper} characterizes the minimal skeleton of any superelliptic curve $y^{n}=f(x)$ in terms of the minimal skeleton of the marked curve $(\mathbb{P}^{1},B)$, where $B$ is the branch locus of the map $(x,y)\mapsto{x}$. Here one assumes $\mathrm{gcd}(n,\mathrm{char}(k))=1$. If $\mathrm{char}(k)\neq{2}$, then this material thus gives the minimal skeleton for curves of genus two. There are seven different tree types for $(\mathbb{P}^{1},B)$ which can be found in Figure \ref{TreesSixLeaves}. For each of these trees, there is a unique degree-$2$ 
covering ramified over its leaves (see \cite{DecompFund1}) and the covering curve has genus $2$. The tree types in Figure \ref{TreesSixLeaves} now match up with the corresponding reduction types given in Figure \ref{Genus2Mod}. 
\begin{rem}

In \cite{RSS}, tropicalizations of various classical moduli spaces are studied. They in particular study the tropicalization of the moduli space of curves of genus two and the moduli space of projective lines with six marked points. 
There is a minor technical issue involving the base field however in \cite[Theorem 5.3]{RSS}. Namely, they assume in the beginning of the paper 
that  $\mathrm{char}(K)=0$, which is not enough to completely remove the problem of wild ramification for the hyperelliptic covering $X\to\mathbb{P}^{1}$ as we can take fields of residue characteristic two, e.g., $\mathbb{C}_{2}$.  In fact, their result becomes false in this case: the curve $y^2=x^5+x$ has tree type IV with respect to our labeling, but its reduction type is V, see \cite[Exemple 2]{liu}. In particular, we find that the commutative diagram given in \cite[Theorem 5.3.(a)]{RSS} does not commute. To repair this issue, one can just assume $\mathrm{char}(k)\neq{2}$ as in this paper. 
\end{rem}

\begin{rem}\label{CuetoMarkwigRemark}
In \cite{CuetoMarkwig2016}, curves of genus two are studied from a tropical point of view. To find the skeleton of a curve of genus two, they use the notion of a faithful tropicalization. In short, this is an embedding of $X$ into a toric variety 
such that its tropicalization defines an isometry of a skeleton of $X$ into the tropicalization of the toric variety. This method does not directly work for two reduction types however, see \cite[Theorem 1.2]{CuetoMarkwig2016}. They also obtain 
a subversion of our Theorem \ref{MainThm1}, determining conditions on the tropical Igusa invariants for some of the reduction types and giving formulas for the edge lengths in some cases. Our results and proofs do not encounter any of the difficulties above: we can determine the skeleton of any curve and we have a complete description of the skeleta (including the edge lengths) in terms of our tropical Igusa invariants. This method does not give a faithful embedding of the skeleton, but these can be produced abstractly using \cite{BPR11}. Theorem \ref{MainThm1} does determine the invariants (e.g. the cycle lengths) of a faithful tropicalization of $X$ however, since these only rely on the abstract minimal skeleton.  
\end{rem}

\subsubsection{$\mathrm{PGL}_{2}$ and trees}\label{SectionTrees}

As before, let $K$ be a complete, algebraically closed non-archimedean field. The automorphism group of the projective line $\mathbb{P}^{1}:=\mathbb{P}^{1}_{K}$ is $\mathrm{PGL}_{2}(K)$ and consequently by GAGA the automorphism group of the Berkovich analytification $\mathbb{P}^{1,\mathrm{an}}$ is also $\mathrm{PGL}_{2}(K)$. In projective coordinates, a matrix $\tau=\begin{pmatrix}
a & b\\
c & d
\end{pmatrix} $ acts by $\tau([x:y])=[ax+by:cx+dy]$. We also refer to these as M\"{o}bius transformations. Throughout this paper, we will often use the fact that for any three distinct points $P_{i}$ in $\mathbb{P}^{1}$, we can find a unique $\tau\in\mathrm{PGL}_{2}(K)$ such that $\tau(P_{1})=0$, $\tau(P_{2})=1$ and $\tau(P_{3})=\infty$. 
This will greatly reduce some of our computations.

Consider the hyperbolic space $\mathbf{H}_{0}:=\mathbf{H}_{0}(\mathbb{P}^{1,\mathrm{an}})$ as in \cite[Section 5]{BPRa1}. There is a natural skeletal metric on $\mathbf{H}_{0}$ and every automorphism of $\mathbb{P}^{1,\mathrm{an}}$ acts as an isometry by the results in \cite[Section 2]{BPRa1}. We now add a set of marked points $D$ to $\mathbb{P}^{1}$. If $|D|\geq{2}$, then the marked curve $(\mathbb{P}^{1},D)$ satisfies $\chi(\mathbb{P}^{1},D)\leq{0}$ and we thus have a unique minimal skeleton $\Sigma_{D}:=\Sigma(\mathbb{P}^{1},D)$. For any automorphism $\tau\in\mathrm{PGL}_{2}(K)$, we then have $\tau(\Sigma_{D})=\Sigma_{\tau(D)}$. Combinatorially, we can just replace every leaf $v\in\Sigma_{D}$ with $\tau(v)$ to find $\tau(\Sigma_{D})$. 

\begin{figure}[h!]
\centering
\includegraphics[height=5cm]{{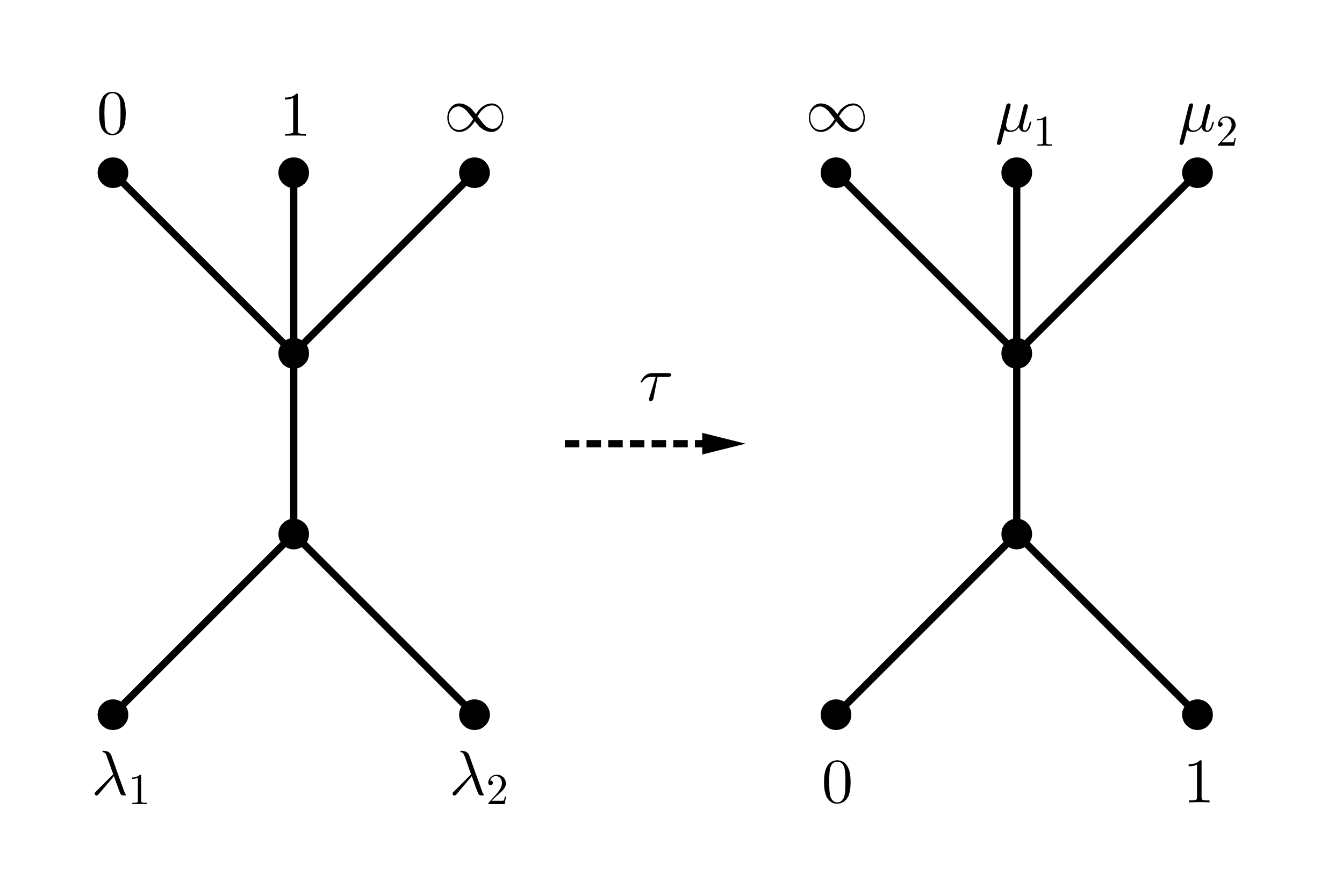}}
\caption{\label{TreeExampleFiveLeaves}
The transformation of minimal skeleta in Example \ref{ExampleFiveLeaves}. 
} 
\end{figure}

\begin{exa}\label{ExampleFiveLeaves}
Consider the type-$1$ points $0,1,\infty,\lambda_{1},\lambda_{2}$ as in Figure \ref{TreeExampleFiveLeaves}. There is then a unique automorphism $\tau$ such that $\tau(\lambda_{1})=0$, $\tau(\lambda_{2})=1$ and $\tau(0)=\infty$ and the corresponding matrix is $\begin{pmatrix}
\lambda_{2} & -\lambda_{1}\lambda_{2} \\
\lambda_{2}-\lambda_{2} & 0
\end{pmatrix}$. The group of three leaves after this transformation is given by $\infty$, $\mu_{1}:=\lambda_{2}(1-\lambda_{1})/(\lambda_{2}-\lambda_{1})$ and $\mu_{2}:=\lambda_{2}/(\lambda_{2}-\lambda_{1})$, see Figure \ref{TreeExampleFiveLeaves}. 
\end{exa}

\subsection{Igusa invariants}\label{SectionIgusa}

We now recall some concepts and results from \cite{Igusa1960} and \cite{liu} on invariants and moduli of curves of genus two.  Let $T=K[u_{0},...,u_{n}]$ and write $f=u_{0}x^{n}+u_{1}x^{n-1}y+...+u_{n}y^{n}\in{T[x,y]}$ for the generic binary form of $(x,y)$-degree $n$ over an algebraically closed field $K$. We define a degree function on $T$ by setting $\mathrm{deg}(u_{i})=i$. We also refer to this as the $T$-degree. As before, any matrix $\sigma=\begin{pmatrix}
a & b \\
c & d
\end{pmatrix}$ in $\mathrm{PGL}_{2}(K)$ acts on $f$ through M\"{o}bius transformations. By comparing the entries of $\sigma(f)$ and $f$, we then also obtain an action on $T$. An {\it{invariant}} of index $g$ is a homogeneous polynomial $H\in{T}$ such that for any matrix $\sigma\in\mathrm{PGL}_{2}(K)$, we have $\sigma(H)=\mathrm{det}(\sigma)^{g}H$. The degree of an invariant is its degree as an element of $T$. In this paper we are interested in invariants of binary forms
$$f=u_{0}x^{6}+u_{1}x^{5}y...+u_{6}y^{6}$$ of $(x,y)$-degree $6$. We moreover restrict ourselves to invariants of even $T$-degree. These invariants form an algebra generated over $\mathbb{Q}$ by elements $A$, $B$, $C$ and $D$ of degrees $2$, $4$, $6$ and $8$ respectively. These are used in \cite{Igusa1960} to define invariants $J_{2i}$ of degree $2i$ for $i=1,2,3,4,5$. We refer to these polynomials as the {\it{Igusa invariants}}.  To calculate them, we use the command $\texttt{JInvariants}(f)$ in the computer algebra package MAGMA.

There is one relation between the $J_{2i}$, namely $J_{4}^2+4J_{8}=J_{2}J_{6}$. We denote the corresponding graded algebra by $\mathbb{Z}[\overline{J}]=\mathbb{Z}[J_{2},J_{4},J_{6},J_{8},J_{10}]$\footnote{This ring is isomorphic to $\mathbb{Z}[X_{2i}]/(X_{4}^2+4X_{8}-X_{2}X_{6})$.} and its projectivization by $\mathcal{X}=\mathrm{Proj}\,\mathbb{Z}[\overline{J}]$. We can connect these invariants to the moduli space of genus two curves as follows. Write $\mathcal{M}_{2}:(\mathrm{Schemes})\to(\mathrm{Sets})$ for the functor that assigns to a scheme $S$ the set of isomorphism classes of smooth proper morphisms $C\to{S}$ whose fibers are geometrically connected curves of genus two. We denote the corresponding coarse moduli space by $M_{2}$. The Igusa invariants allow us to give this coarse moduli space explicitly: it is the open subscheme $D_{+}(J_{10})\subset{\mathcal{X}}$, see \cite[Theorem 2]{Igusa1960}.

\subsubsection{Tropicalizing the Igusa invariants}\label{SectionTropIgusa}

We now consider the algebra $K[\overline{J}]=\mathbb{Z}[\overline{J}]\otimes_{\mathbb{Z}}{K}$ over a complete algebraically closed non-archimedean field $K$. We write $Z$ for the corresponding scheme. Using the grading $\mathrm{deg}(J_{2i})=2i$, we can view this as a closed subscheme of the weighted projective space $\mathbb{P}(2,4,6,8,10)$. Since weighted projective varieties are toric varieties, we can now define the tropicalization of $Z$ with respect to this embedding, see \cite[Chapter 6]{tropicalbook},  \cite{Gubler2013}, \cite{Payne2009} or \cite{Rab2012}. Unfortunately, this tropicalization does not determine the abstract tropicalization of ${M}_{2}$. 
\begin{exa}
Let $\mathrm{char}(k)\neq{2}$ and consider the three polynomials 
\begin{align*}
f_{a}&=x(x-ta_{1})(x-1)(x-(1+ta_{2}))(x+1),\\
f_{b}&=x(x-tb_{1})(x-1)(x-(1+tb_{2}))(x+1),\\
g&=x(x-1)(x-t^2)(x-(t^2+t^4))(x-(t^2+2t^4)),
\end{align*}
where $v(a_{1})=v(a_{2})=3$, $v(b_{1})=2$ and $v(b_{2})=4$. The valuations of the invariants $J_{2i}$ are the same for $f_{a}$ and $f_{b}$, but the lengths of the edges in the minimal skeleta of the curves $y^2=f_{a}$ and $y^2=f_{b}$ are different. The valuations of the Igusa invariants of $g$ differ from the other two by the vector $(2,4,6,8,10)$. That is, they are equivalent in the corresponding tropical weighted projective space (see Definition \ref{WeightedTropProj}).\footnote{We can also calculate the invariants of the polynomial $g_{1}=1/t^{8}g(t^2x)$. The valuations of the Igusa invariants of $g_{1}$ are then the same as the valuations of the Igusa invariants of $f_{a}$ and $f_{b}$.} Note that the reduction type for $f_{a}$ and $f_{b}$ is III, but the reduction type for $g$ is VI. We thus see that the valuations of these Igusa invariants can in general detect neither the edge lengths, nor the reduction type.    
 
\end{exa}   

To obtain the abstract tropicalization of $M_{2}$, we extend our set of Igusa invariants as in \cite{liu}. Define
\begin{eqnarray*}
I_{2}&:=&12^{-1}J_{2},\\
I_{4}&:=&J_{2}^{2}-2^3\cdot{3}\cdot{}J_{4},\\
I_{6}&:=&J_{6},\\
I_{8}&:=&J_{8},\\
I_{12}&:=&-2^{3}J_{4}^{3}+3^{2}J_{2}J_{4}J_{6}-3^{3}J_{6}^{2}-J_{2}^{2}J_{8}.
\end{eqnarray*} 
If we remove the redundant factors, this gives an embedding of $\mathrm{Proj}(K[\overline{J}])$ into $\mathbb{P}(2,2,4,4,6,8,10,12)$. To express this tropicalization more neatly, we introduce the corresponding weighted tropical projective space.  

\begin{mydef}\label{WeightedTropProj}
Let $\overline{\mathbb{R}}=\mathbb{R}\cup\{\infty\}$ and consider the punctured tropical plane $(\mathbb{T}\mathbb{A}^{n})^{*}=\overline{\mathbb{R}}^{n}\backslash{\{(\infty,...,\infty)\}}$. Let $(a_{i})\in\mathbb{N}^{n}$ be fixed. We define an action of $\mathbb{R}$ on the punctured tropical plane by 
\begin{equation}
\lambda\odot(x_{i}):=(x_{i})+\lambda\cdot{}(a_{1},...,a_{n}).
\end{equation} 
The tropical weighted projective space $\mathbb{T}\mathbb{P}(a_{1},...,a_{n})$ is the set-theoretic quotient of the tropical punctured plane by this action. 

\end{mydef}

Using this definition with $(a_{i})=(2,2,4,4,6,8,10,12)$, we arrive at our definition of the tropical Igusa invariants.

\begin{mydef}\label{DefinitionTropicalIgusaInvariants}
{\bf{[Tropical Igusa Invariants]}}
The valuations of the invariants $J_{2i}$ and $I_{2i}$ 
are the tropical Igusa invariants. We view these as functions $Z({K})\to\mathbb{T}\mathbb{P}(2,2,4,4,6,8,10,12)$. 
\end{mydef}       

\begin{rem}
If we are not interested in non-archimedean fields of residue characteristic $\mathrm{char}(k)=2,3$, then we do not need to consider $I_{2}$. The corresponding embedding of $\mathrm{Proj}(K[\overline{J}])$ then lands in $\mathbb{P}(2,4,4,6,8,10,12)$ and the target space for the tropicalization becomes $\mathbb{T}\mathbb{P}(2,4,4,6,8,10,12)$. 
\end{rem}

\begin{rem}
We go back one step and review a similar construction for binary quartics and elliptic curves. As expected, this gives the tropical $j$-invariants studied in \cite{katz_markwig_markwig_2009} and \cite{FaithTropEll}. For simplicity, we assume that $\mathrm{char}(k)\neq{2,3}$ here.  
The invariant ring for binary forms $f=c_{0}x^{4}+...+c_{4}y^{4}$ of $(x,y)$-degree four admit two natural generators which are usually denoted by $c_{4}$ and $c_{6}$. The discriminant $\Delta$ of $f$ can then also be expressed in terms of these generators. For instance, in terms of the formulae given in \cite{silv1} we have $1728\Delta=c_{4}^3-c_{6}^2$. Here we view $0\cdot{}x^{4}+4x^3y+b_{2}x^2y^2+2b_{4}xy^3+b_{6}y^4$ as a degenerate quartic with one zero at infinity (this is generic enough for our purposes). Consider the scheme $Z=\mathrm{Proj}\,{K[c_{4},c_{6},\Delta]}$, where $c_{4}$, $c_{6}$ and $\Delta$ have degrees $4$, $6$ and $12$ respectively. 
We view this as a closed subscheme of the weighted projective space $\mathbb{P}(4,6,12)$. The tropicalization of $\mathbb{P}(4,6,12)$ is $\mathbb{T}\mathbb{P}(4,6,12)$ and the functions $\mathrm{val}(c_{4})$, $\mathrm{val}(c_{6})$ and $\mathrm{val}(\Delta)$ are the tropical invariants for the quartic. The tropical $j$-invariant is then the rational function $3\mathrm{val}(c_{4})-\mathrm{val}(\Delta)$ defined on the open subscheme $D_{+}(\Delta)$ of $Z$. 
\end{rem}

\begin{rem}
Returning to the Igusa invariants, we now explain where $I_{4}$ and $I_{12}$ come from. If we evaluate these invariants at the degenerate sextic $f=0\cdot{}x^6+0\cdot{}x^5y+x^4y^2+...$, we obtain the $c_{4}$ and $\Delta$ invariants for $f$ considered as a quartic, see the paragraph before \cite[D\'{e}finition 1]{liu}. We will not explicitly need this fact in our proof of Theorem \ref{MainThm1}. 
\end{rem}

We now define the subspaces of $\mathbb{T}(2,2,4,4,6,8,10,12)$ that determine the various reduction types of genus two curves. These subspaces will be defined by homogeneous linear forms in the tropical Igusa invariants. To give linear forms that work in all residue characteristics, we first define the following: 
\begin{eqnarray*}
\epsilon:=\begin{cases}
1\text{ if char}(k)\neq{2,3},\\
3\text{ if char}(k)=3,\\
4\text{ if char}(k)=2.
\end{cases}
\end{eqnarray*}  
The linear forms are now given by 
\begin{eqnarray*}
w_{1,i}&=&5v(J_{2i})-iv(J_{10}),\\
w_{2,i}&=&6v(J_{2i})-iv(I_{12}),\\
w_{3,i}&=&2v(J_{2i})-iv(I_{4}),\\
w_{3,x}&=&v(I_{12})-3v(I_{4}),\\
w_{3,y,1}&=&v(J_{4})-v(I_{4}),\\
w_{3,y,2}&=&2v(J_{6})-3v(I_{4}),\\
w_{4,i}&=&w_{3,i},\\
w_{5,1}&=&3\epsilon{v(I_{4})}-\epsilon{}v(J_{10})-v(I_{2\epsilon}),\\
w_{5,2}&=&\epsilon{}v(I_{12})-\epsilon{}v(J_{10})-v(I_{2\epsilon}),\\
\end{eqnarray*}
\begin{eqnarray*}
w_{6,1}&=&3v(I_{4})-v(I_{12}),\\
w_{6,2}&=&\epsilon{}v(J_{10})+v(I_{2\epsilon})-\epsilon{}v(I_{12}),\\
w_{7,1}&=&-w_{6,1},\\
w_{7,2}&=&\epsilon{}v(J_{10})+v(I_{2\epsilon})-3\epsilon{}v(I_{4}),\\
w_{2c,1}&=&\epsilon{v}(I_{4})-2v(I_{2\epsilon}),\\
w_{2c,2}&=&\epsilon{v}(J_{10})-5v(I_{2\epsilon}),\\
w_{2c,3}&=&\epsilon{v}(I_{12})-6v(I_{2\epsilon}).\\
\end{eqnarray*}
Here the first subscript for all but the last three linear forms corresponds to the reduction type. 
Note that all these functions are tropically homogeneous, so that they give rational functions on open subsets of the tropical weighted projective space $\mathbb{T}\mathbb{P}(2,2,4,4,6,8,10,12)$. It will follow from the proof of Theorem \ref{ProofMainThm2} that these are well defined. For instance, for curves of reduction type II, we have that $I_{12}$ is nonzero, so $v(I_{12})\neq{\infty}$ and we can thus evaluate the $w_{2,i}$. If a given set of Igusa invariants has $I_{12}=0$, then we write $w_{2,i}<0$, so that the conditions for II in Theorem \ref{MainThm2} are not satisfied.

\begin{theorem}
\label{MainThm2}
{\bf{[Main Theorem I]}}
The reduction types for curves of genus two are described by the following conditions on the tropical Igusa invariants.  
\begin{itemize}
\item[I.] The reduction type is I if and only if $w_{1,i}\geq{0}$ for every $1\leq{i}\leq{5}$. 
\item[II.] The reduction type is II if and only if $w_{2,i}\geq{0}$ for every $1\leq{}i\leq{4}$ and $w_{2,5}>{0}$.
\item[III.] The reduction type is III if and only if  $w_{3,i}\geq{0}$ for every $i\leq{5}$, $w_{3,5}>0$, $w_{3,x}\geq{0}$ and either $w_{3,y,1}=0$ or $w_{3,y,2}=0$.
\item[IV.] The reduction type is IV if and only if $w_{4,i}>0$ for every $2\leq{i}\leq{5}$.
\item[V.] The reduction type is V if and only if $w_{2c,i}>0$ for $i=1,2,3$ and $w_{5,1},w_{5,2}\geq{0}$. 
\item[VI.] The reduction type is VI if and only if $w_{2c,i}>0$ for $i=1,2,3$, $w_{6,1}\geq{0}$ and $w_{6,2}>0$.
\item[VII.] The reduction type is VII if and only if $w_{2c,i}>0$ for $i=1,2,3$ and $w_{7,1},w_{7,2}>0$. 
\end{itemize}
\end{theorem}

The proof of this theorem will be given in Section \ref{ProofMainThm2}. 

\begin{rem}
We will only prove this when $\mathrm{char}(k)\neq{2}$. \cite[Th\'{e}or\`{e}me 1]{liu} states 
 that the same conditions also work for fields of residue characteristic two,
but we haven't verified this ourselves. The definition for $\epsilon$ when $\mathrm{char}(k)=2,3$ is an artifact of certain anomalous congruences on Igusa invariants, see part (B) of the proof of \cite[Th\'{e}or\`{e}me 1]{liu}. The reader who is not interested in these cases can just assume $\epsilon=1$ throughout the rest of the paper.  

\end{rem}

\begin{rem}
The discretely valued case proven in \cite{liu} can be recovered from the algebraically closed case we prove here. Indeed, if $X$ has stable reduction over $K$, then the skeleton is stable under taking the base change to the completion of the algebraic closure.   
\end{rem}

\begin{rem}
Th\'{e}or\`{e}me 1 in \cite{liu} also determines the $j$-invariants of the various elliptic curves in the reduction types in terms of the Igusa invariants. For the sake of brevity, we have omitted these formulas but they can be easily checked using the universal families given in Section \ref{ProofMainThm2}.   
\end{rem}

The tropical Igusa invariants also determine the edge lengths of the minimal skeleton of a curve of genus two. This is the content of the second part of our main theorem. 

\begin{theorem}\label{MainThm3}
{\bf{[Main Theorem II]}}
The edge lengths of the minimal skeleton of a curve $X$ of genus two are given by the following. 
\begin{itemize}
\item[I.] There are no non-trivial edges. 
\item[II.]  There is one non-trivial edge of length 
$\ell=w_{2,5}/6$.

\item[III.] 
Let $\ell_{1}\leq{\ell_{2}}$ be the lengths of the two non-trivial edges. Then
\begin{eqnarray*}
\ell_{1}&=&\inf\{w_{3,x},w_{3,5}/4\},\\
\ell_{2}&=&w_{3,5}/2-\ell_{1}.
\end{eqnarray*}
\item[IV.] Let $\ell_{1}\leq{\ell_{2}}\leq{\ell_{3}}$ be the lengths of the three non-trivial edges. Let $L=v(J_{10})-5v(J_{2})$, $N=v(I_{12})-6v(J_{12})$ and $M=v(J_{4})-2v(J_{2})$. Then 
\begin{eqnarray*}
\ell_{1}&=&\inf\{L/3,N/2,M\}, \\
\ell_{2}&=&\inf\{(L-\ell_{1})/2,N-\ell_{1}\}, \\
\ell_{3}&=&L-\ell_{1}-\ell_{2}. 
\end{eqnarray*}

\item[V.] 
There is one non-trivial edge of length 
\begin{equation*}
\ell=\dfrac{1}{12\epsilon}\epsilon{}v(J_{10})-5v(I_{2\epsilon}).
\end{equation*}
\item[VI.] 
Let $\ell_{0}$ be the length of the connecting edge and let 
$\ell_{1}$ the length of the cycle. 
Then 
\begin{eqnarray*}
\ell_{0}&=&\dfrac{1}{12\epsilon}(\epsilon{}v(I_{12})-6v(I_{2\epsilon})),\\
\ell_{1}&=&\dfrac{1}{\epsilon}(\epsilon{}v(J_{10})+v(I_{2\epsilon})-\epsilon{}v(I_{12})).
\end{eqnarray*}
\item[VII.] Let $\ell_{0}$ be the length of the connecting edge and let $\ell_{1}\leq{\ell_{2}}$ be the lengths of the cycles. 
Then
\begin{eqnarray*}
\ell_{0}&=&\dfrac{1}{4\epsilon}(\epsilon{}v(I_{4})-2v(I_{2\epsilon})),\\
\ell_{1}&=&\inf\{w_{7,1},\dfrac{1}{2\epsilon}w_{7,2}\},\\
\ell_{2}&=&\dfrac{1}{\epsilon}(w_{7,2})-\ell_{1}.
\end{eqnarray*}

\end{itemize}
\end{theorem}

The proof of this theorem will be given in Section \ref{ProofMainThm3}.

\begin{rem}\label{RemarkComponentGroup}
Suppose that $X/K$ is a curve of genus two over a discretely valued complete non-archimedean field $K$. If $X$ has stable reduction over $K$, then the formulas above can be used to determine the component group of the N\'{e}ron model $\mathcal{J}$ of the Jacobian $J$ of $X$, see \cite{Baker2008}. We refer the reader to \cite{liu} for explicit formulas.   
\end{rem}

We now discuss how this material connects to the abstract tropicalizations studied in \cite{ACP2015}. On the algebraic side, there is the coarse moduli space $M_{2}$, see the beginning of Section \ref{SectionIgusa}. On the tropical side we have 
$M^{\mathrm{trop}}_{2}$, which is the moduli space of tropical curves of genus two. We quickly recall its definition here, see \cite{ACP2015} for more details.  Let $\mathbf{G}$ be a weighted stable graph of genus two. We can vary the lengths of the edges in $\mathbf{G}$ to obtain the local moduli space $\sigma^{o}_{\mathbf{G}}=\mathbb{R}_{>0}^{|E(\mathbf{G})|}$ for tropical curves of type $\mathbf{G}$. We write $\sigma_{\mathbf{G}}$ for the closed cone corresponding to $\sigma^{o}_{\mathbf{G}}$. The boundary  $\sigma_{\mathbf{G}}\backslash{\sigma^{o}_{\mathbf{G}}}$ also has a moduli space interpretation: it corresponds to reduction types $\mathbf{G}'$ obtained by contracting edges or cycles $e\in{E(\mathbf{G})}$ with $\ell(e)=0$. For every contraction $\varpi:\mathbf{G}\to\mathbf{G}'$, we then obtain an inclusion of closed cones $j_{\varpi}:\sigma_{\mathbf{G}'}\to\sigma_{\mathbf{G}}$. Here we also allow isomorphisms of graphs, which are in a sense degenerate contractions.  
The tropical moduli space for curves of genus two is then defined as the colimit $M_{2}^{\mathrm{trop}}=\lim(\sigma_{\mathbf{G}},j_{\varpi})$ over all weighted contractions. Locally this limit is quite concrete: 
the image of $\sigma^{0}_{\mathbf{G}}$ in $M_{2}^{\mathrm{trop}}$ is 
$\sigma^{o}_{\mathbf{G}}/\mathrm{Aut}(\mathbf{G})$. Other examples of contraction morphisms for tropical curves of genus two can be found in Figure \ref{Genus2Mod}. 
\begin{exa}
Let $\mathbf{G}$ be the graph corresponding to reduction type IV. The corresponding open moduli space is $\sigma^{o}_{\mathbf{G}}=\mathbb{R}_{>0}^{3}$ and its image in $M_{2}^{\mathrm{trop}}$ is $\sigma^{o}_{\mathbf{G}}/\mathrm{Aut}(\mathbf{G})=\mathbb{R}_{>0}^{3}/S_{3}$. Consider the sublocus of $\sigma^{o}_{\mathbf{G}}$ outside the diagonals, i.e., consider $(x,y,z)\in\mathbb{R}_{>0}^{3}$ with $x$, $y$ and $z$ distinct. Its image in the quotient can then be uniquely represented by triples $(x,y,z)$ with $x<y<z$. These kinds of decompositions will also play a role in the proof of Theorem \ref{MainThm1}, see \ref{CaseIVEdges}.     
\end{exa}

We can connect the algebraic approach to the tropical approach by introducing an abstract tropicalization map. We start with a point $P\in{M_{2}(K)}$, which corresponds to a smooth curve $X/K$ of genus $2$. The minimal skeleton $\Sigma(X)$ of $X$ is automatically a tropical curve and the tropicalization map $\mathrm{trop}:M_{2}\to{M_{2}^{\mathrm{trop}}}$ sends $X/K$ to $\Sigma(X)$. We now use Theorem \ref{MainThm1} to relate
this abstract tropicalization map to the embedded tropicalization map in \ref{DefinitionTropicalIgusaInvariants}. As mentioned before, the tropical Igusa invariants give us explicit coordinates on the coarse moduli space $M_{2}$, which 
gives an embedding of $M_{2}$ into $\mathbb{P}(2,2,4,4,6,8,10,12)$. We write $\phi$ for this embedding and $\phi^{\mathrm{trop}}$ for its tropicalization. Theorem \ref{MainThm1} then gives a
stratification of $\phi^{\mathrm{trop}}(M_{2})$ in terms of half-spaces. On the other hand, 
$M_{2}^{\mathrm{trop}}$ also has a natural stratification in terms of the reduction types. 
Our findings in Theorem \ref{MainThm1} can now be summarized as follows. 
\begin{cor}\label{CorModuliSpace}
There is a surjective map $i:\phi^{\mathrm{trop}}(M_{2})\to{M_{2}^{\mathrm{trop}}}$ such that $\mathrm{trop}=i\circ\phi^{\mathrm{trop}}$. 
This map preserves the stratifications on both sides.  

\end{cor}
\begin{proof}
We define the map $i$ as follows. Let $P=(v(J_{2i}),v(I_{2i}))$ be a given set of tropical Igusa invariants obtained from a point in $M_{2}(L)$ for a non-archimedean field $L$. By Theorem \ref{MainThm2}, this uniquely determines a reduction type. 
Theorem \ref{MainThm3} then gives us the edge lengths of the corresponding skeleton, which allows us to reconstruct the corresponding tropical curve.  We set $i(P)$ to be this tropical curve. This map obviously commutes with the abstract tropicalization map. 
For the surjectivity, we use $L$-valued points of $M_{2}$, where $K\subset{}L$ is a complete algebraically closed non-archimedean field with value group $\mathbb{R}$. Let $\Sigma$ be a tropical curve of genus $2$. Using the universal families in  
the proof of Theorem \ref{MainThm2} (see Section \ref{ProofMainThm2}), we easily find a curve $X/L$ that tropicalizes to $\Sigma$. 
The image of the corresponding point $P\in{M_{2}(L)}$ in $\phi^{\mathrm{trop}}(M_{2})$ then gives the desired set of tropical Igusa invariants. 
\end{proof}
 
\begin{rem}
The map $i$ not injective. For instance, the fiber over the locus of good reduction is positive dimensional: for any $t\in{K}$ with $v(t)>0$, we set $f=x^5-tx-1$ and find $v(J_{10})=0$ and $v(J_{2})=v(t)$ (here we assume $\mathrm{char}(k)\neq{2,5}$). 
 
\end{rem}

\subsection{Proof of the main theorem}\label{ProofMainThm2}

We will assume throughout this section that $\mathrm{char}(k)\neq{2}$. The main idea in the proofs of theorems \ref{ProofMainThm2} and \ref{ProofMainThm3} is as follows. As we saw earlier, there are exactly seven reduction types for curves of genus two. For every reduction type, we give a (naive) universal family, in the sense that every curve of a given reduction type is isomorphic to a member of this family.  We then check the conditions in Theorem \ref{MainThm2} for the members in this family to obtain the first implication. To obtain the second implication, we show that the conditions are exclusive, in the sense that the conditions for two distinct types cannot hold at the same time. Our proof for Theorem \ref{MainThm3} is similar: we compare the edge lengths in the universal family to the functions in \ref{MainThm3} and see that they coincide. 

\begin{rem}
Throughout our computations below, we will additionally assume $\mathrm{char}(k)\neq{3}$ for simplicity and brevity. The main difference for $\mathrm{char}(k)=3$ is that one has to use a new function that is nonzero on the universal families. More specifically, one uses $I_{6}$ instead of $I_{2}$ for the last three reduction types. This also changes the weight, so the powers for the other invariants need to be changed as well, see the equations for the half-spaces given earlier. 
We leave the calculation of the invariants in this case 
to the reader.

\end{rem}

\subsubsection{Case I: Good reduction}\label{CaseI}

\begin{figure}[h!]
\centering
\includegraphics[height=4.5cm]{{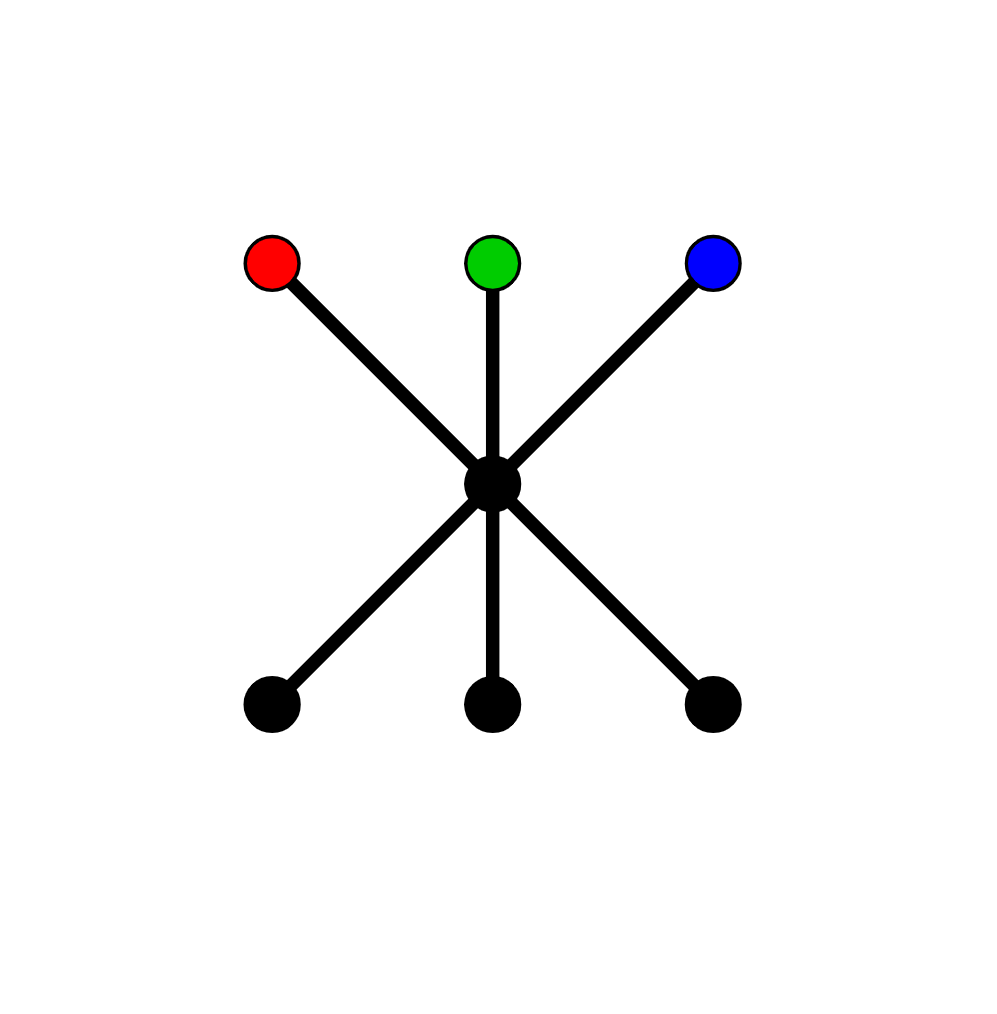}}
\caption{\label{TreeCaseI}
The marked tree used in \ref{CaseI}. The red vertex is $P_{1}$, the green vertex is $P_{2}$ and the blue vertex is $P_{3}$.
} 
\end{figure}

A curve of genus two $X/K$ has good reduction if and only if the phylogenetic type of the tree of the branch locus $B$ of the map $X\to\mathbb{P}^{1}$ is trivial, see \cite[Proposition 3.13]{InvariantsSuper}. We now choose three points $P_{i}$ in $B$ and find a M\"{o}bius transformation $\tau$ such that $\tau(P_{1})=0$, $\tau(P_{2})=1$ and $\tau(P_{3})=\infty$, see Figure \ref{TreeCaseI}. 

 Throughout our proofs, we will use a red vertex for $P_{1}$, a green vertex for $P_{2}$ and a blue vertex for $P_{3}$. Our choice of $\tau$ now implies that the other branch points are mapped to points $\lambda_{i}\in{K}$ 
 of valuation zero. We then similarly have  
$v(\lambda_{i}-1)=0$ and $v(\lambda_{i}-\lambda_{j})=0$. 
 We thus see that $X$ can be represented as 
\begin{equation}
y^2=x(x-1)(x-\lambda_{1})(x-\lambda_{2})(x-{\lambda_{3}}),
\end{equation} 
where the $\lambda_{i}$ satisfy the above conditions. 
To put this in analytic terms, consider the scheme $N$ defined by the ring $K[\lambda_{i}^{\pm},(1-\lambda_{i})^{\pm},(\lambda_{i}-\lambda_{j})^{\pm}]$. The extra tropical conditions on the $\lambda_{i}$ define an affinoid subdomain $U$ of $N^{\mathrm{an}}$ and this is our naive moduli space for curves of good reduction. Using our earlier arguments we now see that every curve $X/K$ of good reduction arises from a type-$1$ point of $U$. Conversely, it is not too hard to see that every type-$1$ point of this space corresponds to a curve $X/K$ of genus two with good reduction. We now calculate the tropical Igusa invariants for this universal family and find that $v(J_{10})=0$ (note that $J_{10}$ is just the discriminant of $f(x)$). This directly shows that the inequalities in Theorem \ref{MainThm2} are satisfied. 

\subsubsection{Case II}\label{CaseII}

\begin{figure}[h!]
\centering
\includegraphics[height=4.5cm]{{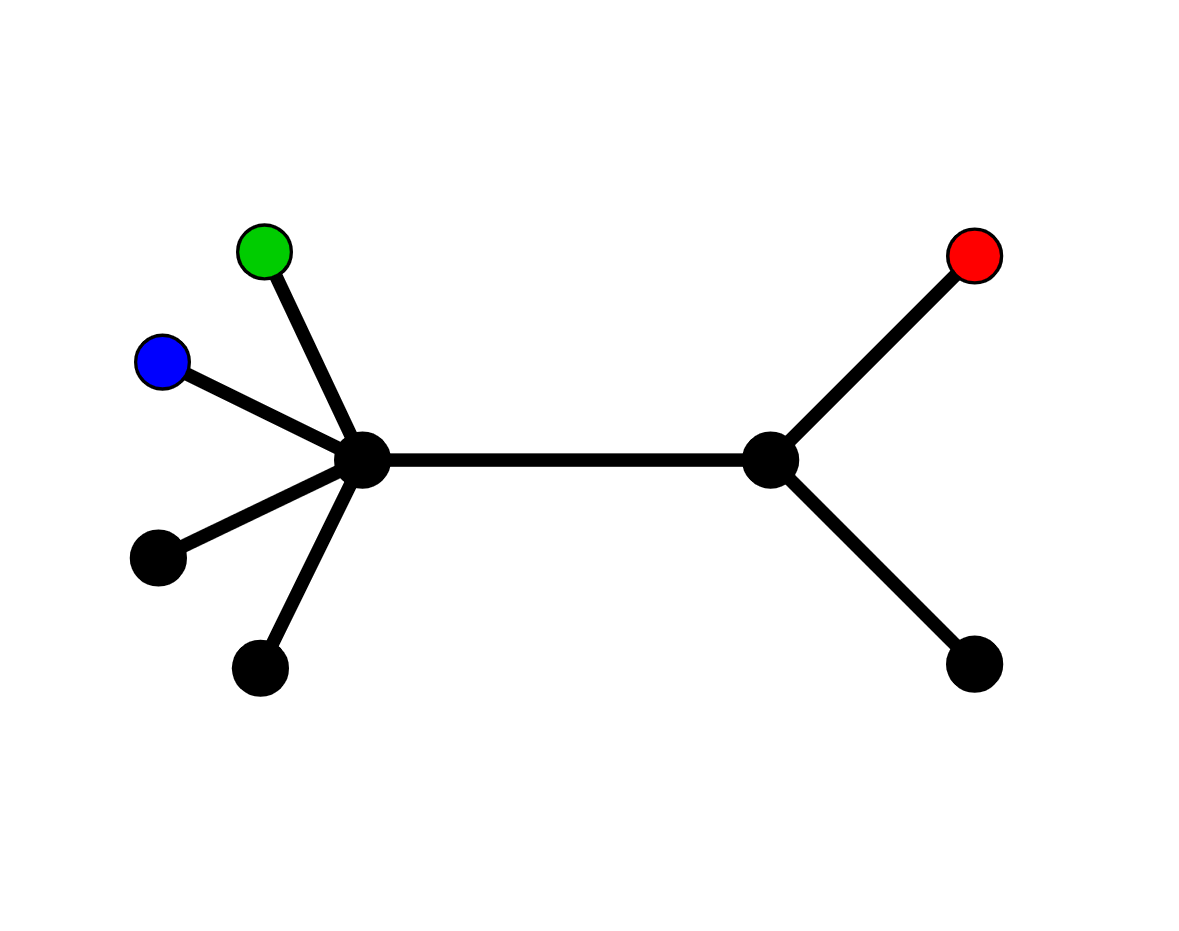}}
\caption{\label{TreeCaseII}
The marked tree used in \ref{CaseII}. The red vertex is $P_{1}$, the green vertex is $P_{2}$ and the blue vertex is $P_{3}$. 
} 
\end{figure}

We choose $\{P_{1},P_{2},P_{3}\}$ as in Figure \ref{TreeCaseII} and   
find a M\"{o}bius transformation $\tau$ such that $\tau(P_{1})=0$, $\tau(P_{2})=1$ and $\tau(P_{3})=\infty$.

 Our universal family is then given by the equation $y^2=f(x)$, where 
\begin{equation}
f(x)=x(x-1)(x-\lambda_{1}{t})(x-\lambda_{2})(x-\lambda_{3}). 
\end{equation}
Here we evaluate $t$ at $\varpi^{k}$ for $k>0$. The equations for the $\lambda_{i}$ are $v(\lambda_{i})=0$ for every $i$, and $v(1-\lambda_{i})=0=v(\lambda_{i}-\lambda_{j})$ for $i=2,3$. We now calculate the reduction modulo $t$ of the Igusa invariant $I_{12}$. This factorizes as 
\begin{equation}
\overline{I}_{12}=1/256\cdot{}\lambda_{3}^6(\lambda_{3}-1)^2\lambda_{2}^6(\lambda_{2}-1)^2(\lambda_{2}-\lambda_{3})^2. 
\end{equation} 
By the assumptions imposed on our $\lambda_{i}$, we find that $\overline{I}_{12}\neq{0}$. This directly implies that the conditions in the theorem are satisfied.

\subsubsection{Case III}\label{CaseIII}

\begin{figure}[h!]
\centering
\includegraphics[height=4cm]{{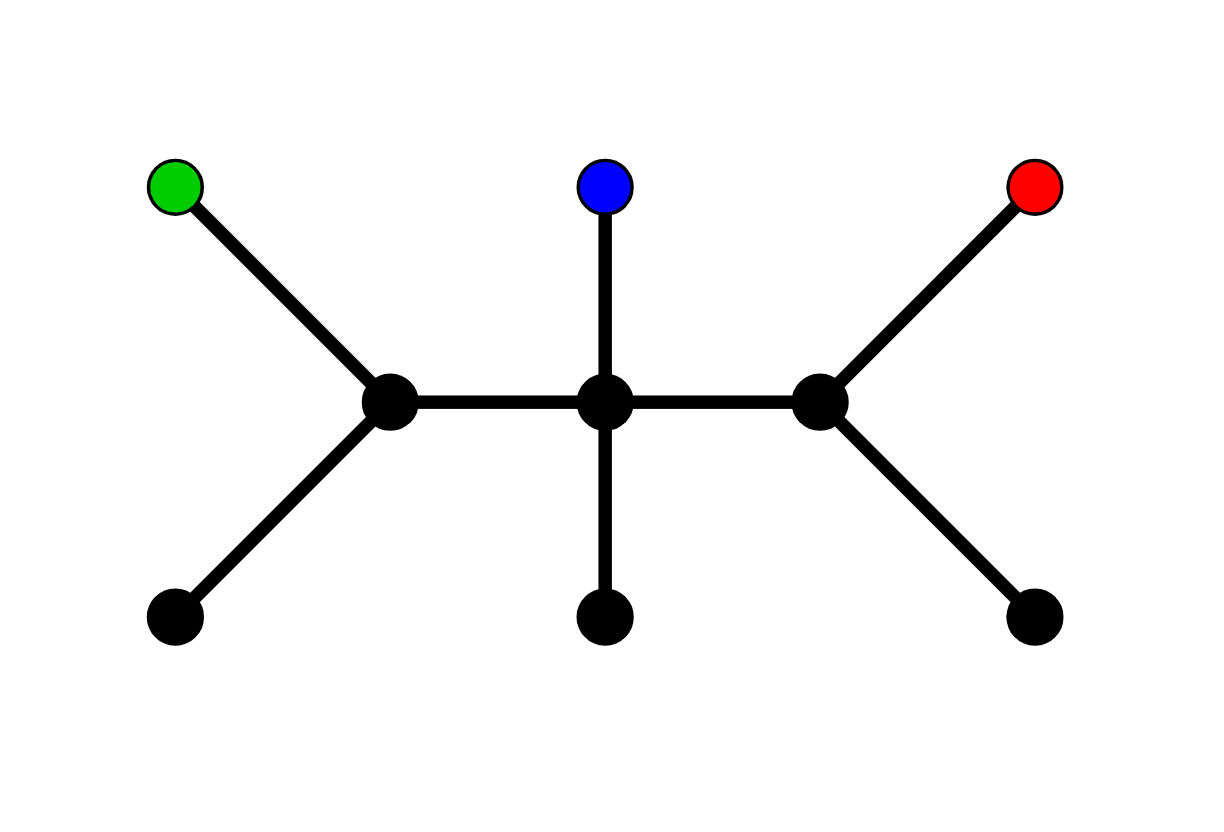}}
\caption{\label{TreeCaseIII}
The marked tree used in \ref{CaseIII}. The red vertex is $P_{1}$, the green vertex is $P_{2}$ and the blue vertex is $P_{3}$. 
} 
\end{figure}

We choose $\{P_{1},P_{2},P_{3}\}$ as in Figure \ref{TreeCaseIII} and   
find a M\"{o}bius transformation $\tau$ such that $\tau(P_{1})=0$, $\tau(P_{2})=1$ and $\tau(P_{3})=\infty$. 

This gives the universal family $y^2=f(x)$ with 
\begin{equation}
f(x)=x(x-1)(x-(1+t_{1}\lambda_{1}))(x-\lambda_{2})(x-t_{2}\lambda_{3}).
\end{equation}
Here $t_{1}$ and $t_{2}$ are evaluated at $\varpi^{k_{i}}$ for $k_{i}=v(t_{i})>0$. We calculate the reduction of $I_{4}$ modulo $(t_{1},t_{2})$ and find 
\begin{equation}
\overline{I}_{4}=(\lambda_{2}-1)^2\lambda_{2}^2,
\end{equation}
which is nonzero by our assumptions. We furthermore have
\begin{align*}
\overline{J}_{4}&=1/16\lambda_{2}^2 - 1/16\lambda_{2} + 3/128,\\
\overline{J}_{6}&=1/256\lambda_{2}^2 - 1/256\lambda_{2} + 1/1024.
\end{align*}
The ideal generated by $\overline{J}_{4}$ and $\overline{J}_{6}$ in $k[\lambda_{2}]$ is trivial, so these two forms cannot simultaneously be zero. This gives the conditions in the theorem.  

\subsubsection{Case IV}\label{CaseIV}

\begin{figure}[h!]
\centering
\includegraphics[height=5cm]{{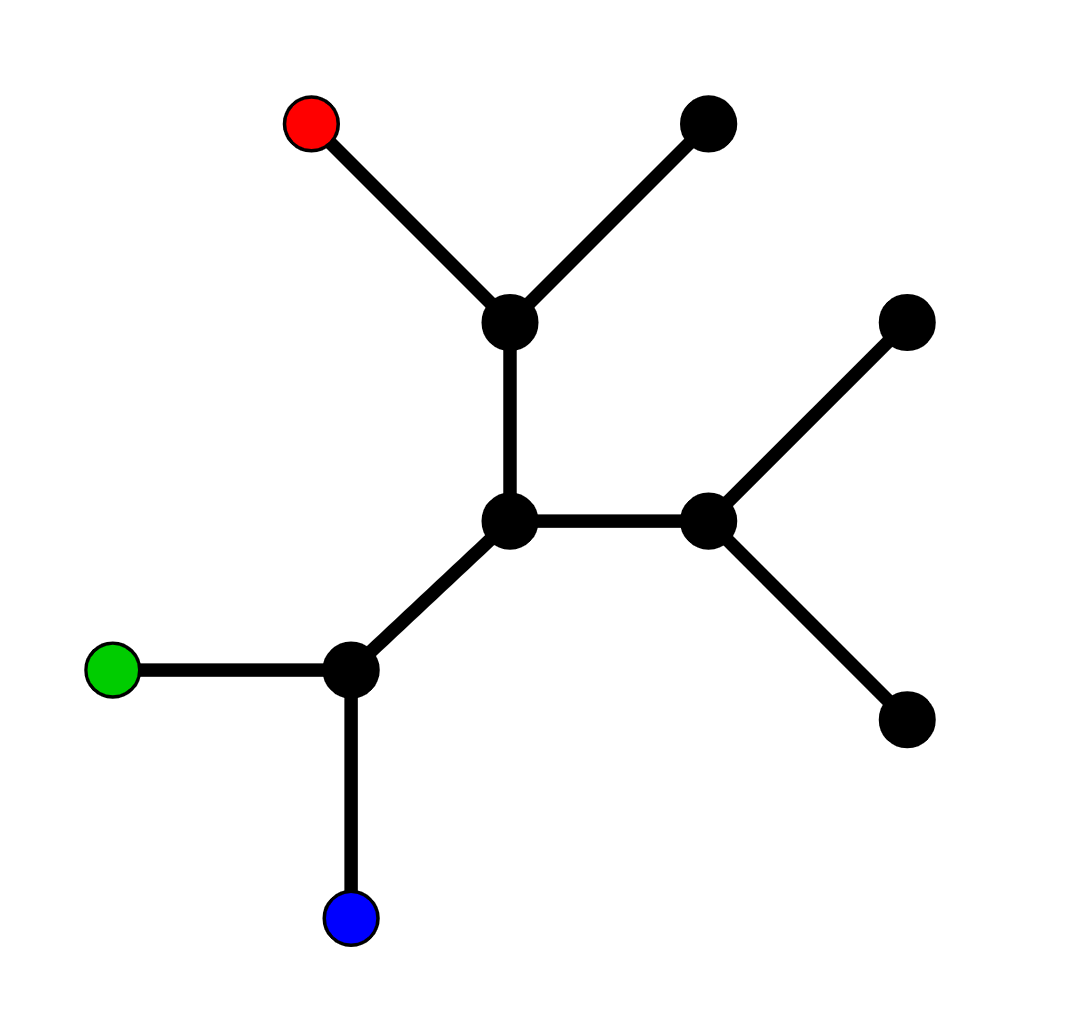}}
\caption{\label{TreeCaseIV}
The marked tree used in \ref{CaseIV}. The red vertex is $P_{1}$, the green vertex is $P_{2}$ and the blue vertex is $P_{3}$. 
} 
\end{figure}

As before, we find a M\"{o}bius transformation that sends the points $P_{1}$, $P_{2}$ and $P_{3}$ in Figure \ref{TreeCaseIV} to $0$, $1$ and $\infty$. 

This gives the universal family $y^2=f(x)$ with 
\begin{align*}
f(x)=x(x-1)(x-h_{1})(x-h_{2})(x-h_{3})
\end{align*}
and 
\begin{align*}
h_{1}&=t_{1}t_{3}\lambda_{1},\\
h_{2}&=t_{1}(\lambda_{2}+\lambda_{3,1}t_{2}),\\
h_{3}&=t_{1}(\lambda_{2}+\lambda_{3,2}t_{2}).
\end{align*}
We then calculate 
 $J_{2}=t_{1}^2m_{2}$, 
$J_{4}=t_{1}^4m_{4}$, $J_{6}=t_{1}^6m_{6}$,
$J_{8}=t_{1}^8m_{8}$,
$J_{10}=t_{1}^{12}t_{2}^2t_{3}^2m_{10}$ and 
$I_{4}=t_{1}^4m^{*}_{4}$. The reduction of $m_{4}^{*}=I_{4}/t_{1}^4$ modulo $(t_{1},t_{2},t_{3})$ is $\lambda_{2}^4$, so $v(I_{4})=4v(t_{1})$. We furthermore find that the reductions of $m_{4}$, $m_{6}$ and $m_{8}$ modulo $(t_{1},t_{2},t_{3})$ are zero. This easily implies the inequalities in the theorem.

\subsubsection{Case V}\label{CaseV}

\begin{figure}[h!]
\centering
\includegraphics[height=4.2cm]{{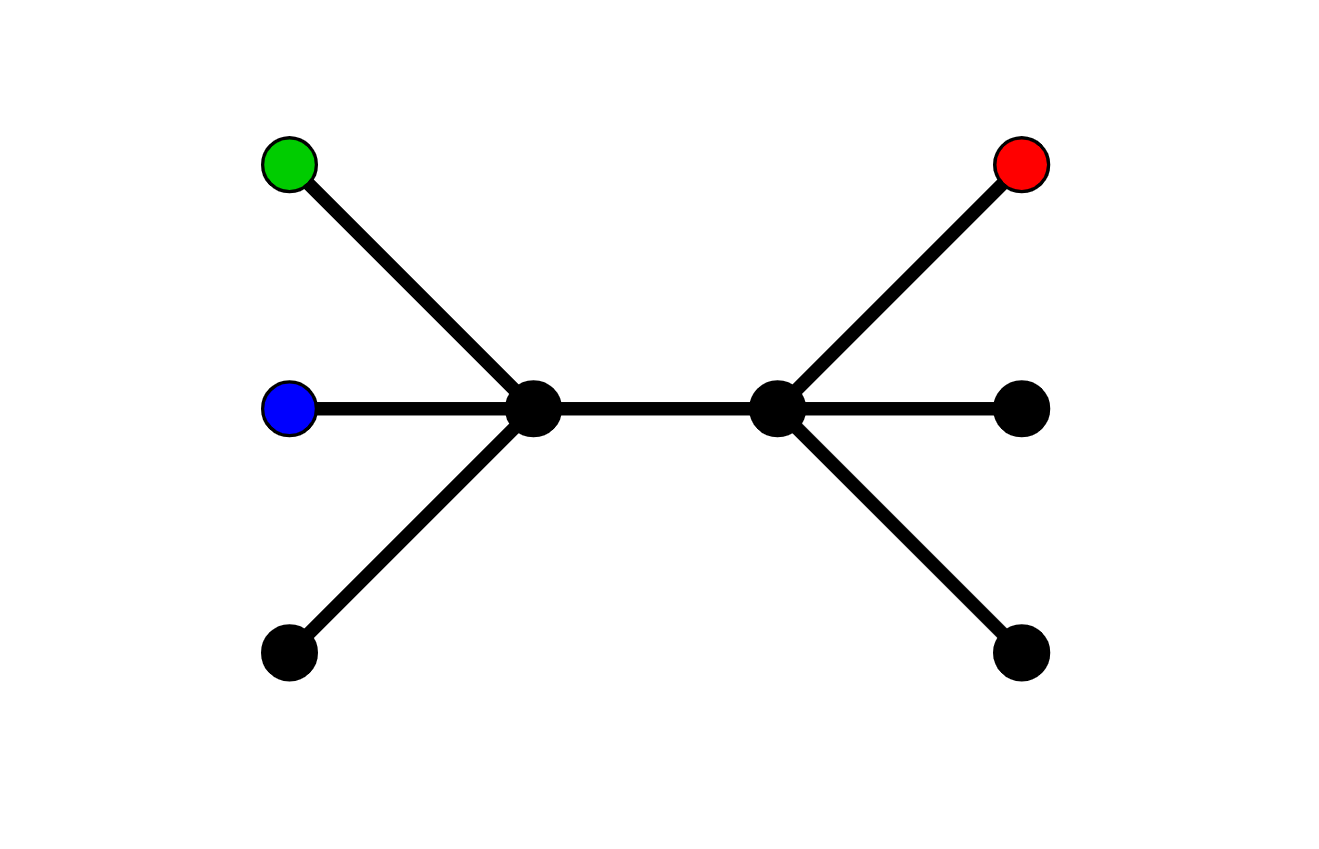}}
\caption{\label{TreeCaseV}
The marked tree used in \ref{CaseV}. The red vertex is $P_{1}$, the green vertex is $P_{2}$ and the blue vertex is $P_{3}$. 
} 
\end{figure}

The universal family derived from Figure \ref{TreeCaseV} is given by $y^2=f(x)$, where 

\begin{align*}
f(x)=x(x-1)(x-h_{1})(x-h_{2})(x-h_{3})
\end{align*}
and 
\begin{align*}
h_{1}&=t_{1}\lambda_{1},\\
h_{2}&=t_{1}\lambda_{2},\\
h_{3}&=\lambda_{3}.
\end{align*}

We first calculate the reduction of $J_{2}$ modulo $t_{1}$, which is $3/4\lambda_{3}^2\neq{0}$. For fields $K$ with $\mathrm{char}{(k)}=3$, this form would be zero, so we would have to replace it by $I_{6}$. At any rate, we now factorize $J_{10}$ and obtain 
\begin{equation}
J_{10}=1/4096\cdot{}t_{1}^{6}\lambda_{3}^2(\lambda_{3}-1)^2\lambda_{2}(\lambda_{2}t_{1}-1)^2(\lambda_{2}t_{1}-\lambda_{3})^2\lambda_{1}^2(\lambda_{1}-\lambda_{2})^2(\lambda_{1}t_{1}-1)^2(\lambda_{1}t_{1}-\lambda_{3})^2.
\end{equation}
We furthermore have $I_{4}=t_{1}^2\cdot{m_{4}}$ and $I_{12}=t_{1}^6\cdot{m_{12}}$ for $m_{4}$ and $m_{12}$ integral. We directly obtain the desired conditions from this.

\subsubsection{Case VI}\label{CaseVI}

\begin{figure}[h!]
\centering
\includegraphics[height=3.5cm]{{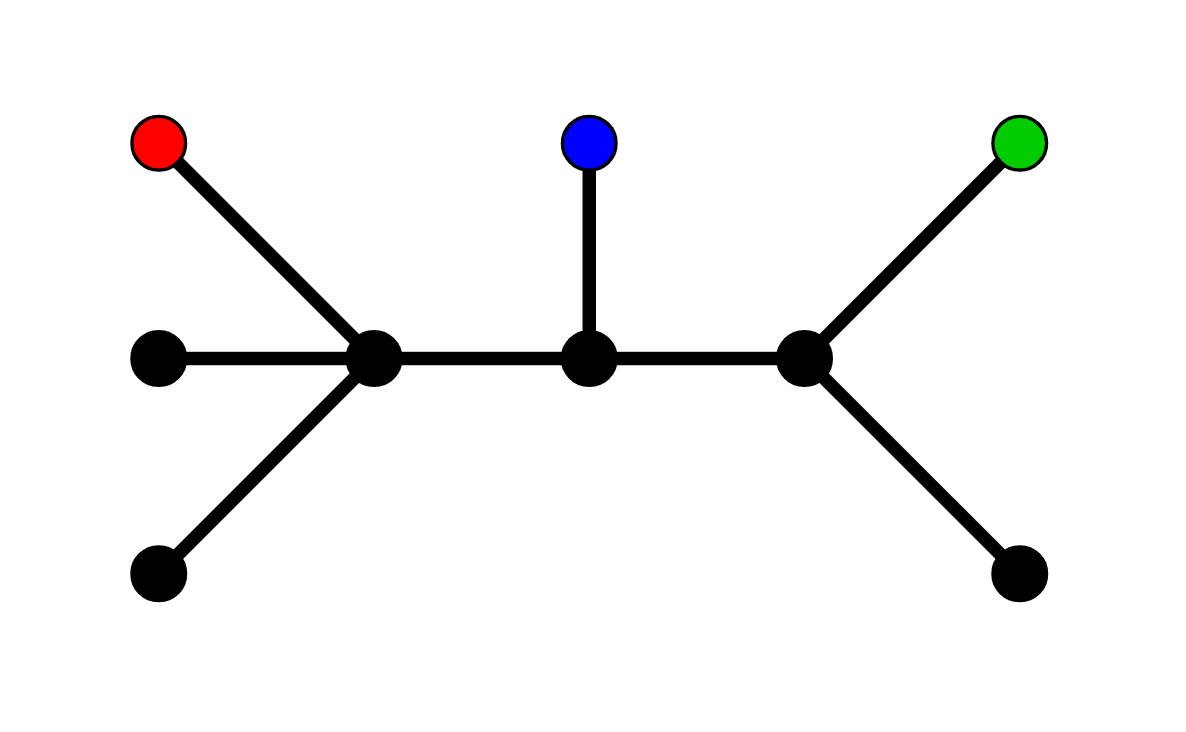}}
\caption{\label{TreeCaseVI}
The marked tree used in \ref{CaseVI}. The red vertex is $P_{1}$, the green vertex is $P_{2}$ and the blue vertex is $P_{3}$. 
} 
\end{figure}

The universal family derived from Figure \ref{TreeCaseVI} is given by $y^2=f(x)$, where 

\begin{align*}
f(x)=x(x-1)(x-h_{1})(x-h_{2})(x-h_{3})
\end{align*}
and 
\begin{align*}
h_{1}&=t_{1}\lambda_{1},\\
h_{2}&=t_{1}\lambda_{2},\\
h_{3}&=1-t_{2}\lambda_{3}.
\end{align*}

The invariant $I_{12}$ is divisible by $t_{1}^{6}$ and the reduction of $I_{12}/t_{1}^6$ modulo $(t_{1},t_{2})$ is $\lambda_{1}^2\lambda_{2}^2(\lambda_{1}-\lambda_{2})^2$, which is not zero. We furthermore have $I_{4}=t_{1}^2m_{4}$ and $J_{10}=t_{1}^6t_{2}^2m_{10}$. The desired inequalities follow directly from this.

  \subsubsection{Case VII}\label{CaseVII}
  
\begin{figure}[h!]
\centering
\includegraphics[height=3.5cm]{{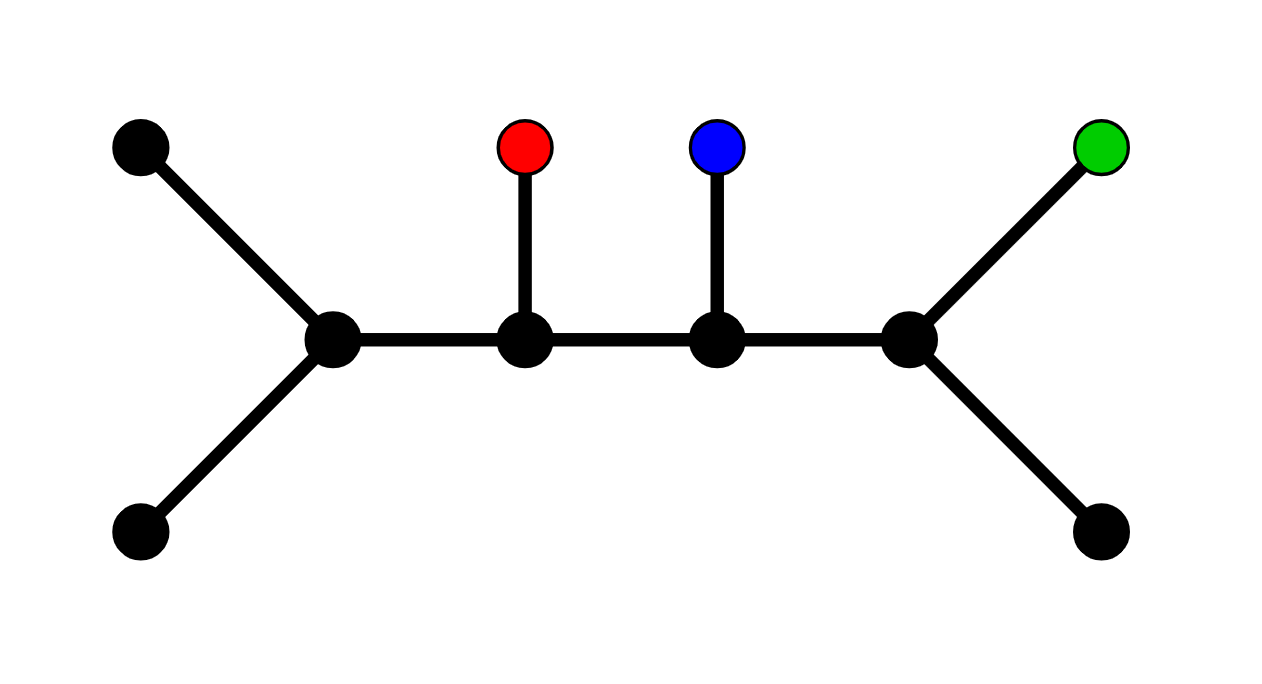}}
\caption{\label{TreeCaseVII}
The marked tree used in \ref{CaseVII}. The red vertex is $P_{1}$, the green vertex is $P_{2}$ and the blue vertex is $P_{3}$. 
} 
\end{figure}

  The universal family derived from Figure \ref{TreeCaseVII} is given by $y^2=f(x)$, where 

\begin{align*}
f(x)=x(x-1)(x-h_{1})(x-h_{2})(x-h_{3})
\end{align*}
and 
\begin{align*}
h_{1}&=t_{1}(\lambda_{1}+\lambda_{2,1}t_{2}),\\
h_{2}&=t_{1}(\lambda_{1}+\lambda_{2,2}t_{2}),\\
h_{3}&=1-t_{3}\lambda_{3}.
\end{align*}

We find that $I_{4}=t_{1}^2m_{4}$, $I_{12}=t_{1}^6m_{12}$ and $J_{10}=t_{1}^6t_{2}^2t_{3}^2m_{10}$. The reduction of $m_{4}$ modulo $(t_{1},t_{2},t_{3})$ is $\lambda_{1}^2$, so it is always nonzero. This gives the desired inequalities. 

\subsubsection{Conclusion}

We now combine the results of the previous sections and prove Theorem \ref{MainThm2}.

\begin{proof}
The previous sections show that the indicated conditions are necessary. For the sufficiency, we show that they define exclusive conditions. We will treat most cases here and leave a few to the reader. Suppose that the conditions for I and II are both satisfied. We can assume that $v(J_{10})=0$. We then have $v(J_{2i})\geq{0}$ and consequently $v(I_{12})\geq{0}$. This contradicts $w_{2,5}>0$, as desired. Suppose that the conditions for I and III are both satisfied. We can assume that $v(J_{10})=0$. The condition for $w_{1,2}$ implies $v(I_{4})\geq{0}$. But then $w_{3,5}\leq{0}$, a contradiction. This also covers the case where I and IV are both satisfied. Suppose that the conditions for II and IV are satisfied. We can suppose that $v(I_{4})=0$, so that $v(J_{2i})>0$ for $2\leq{i}\leq{5}$. This implies 
 $v(J_{2})=0$ by $I_{4}=J_{2}^2-2^3\cdot{3}\cdot{J_{4}}$. The inequality $w_{2,1}\geq{0}$ then gives 
 $v(I_{12})=0$. But then either $v(J_{4})=0$ or $v(J_{6})=0$. This contradicts $w_{4,i}>0$, as desired. We now show that the conditions for III and IV are mutually exclusive. We suppose that $v(I_{4})=0$. Then $v(J_{2i})>0$ for $2\leq{i}\leq{5}$. But this implies that $w_{3,y,1}=0=w_{3,y,2}$, a contradiction.   

Finally, we show that the conditions $w_{2c,i}>0$ for $i=1,2,3$ and the conditions for the first four reduction types are mutually exclusive. We will only do this for $\epsilon=1$ (so $\mathrm{char}(k)\neq{2,3}$), the other cases can be checked in a similar fashion. We start with reduction type I. We can suppose again that $v(J_{10})=0$, so that $v(J_{2i})\geq{0}$. Then $w_{2c,2}$ must be zero, a contradiction. For reduction type II, we suppose that $v(I_{12})=0$, so that $v(J_{2i})\geq{0}$ and $v(J_{10})>0$. But then $w_{2c,3}=0$, a contradiction. For reduction type III, we suppose that $v(I_{4})=0$. We then have $v(J_{2i})\geq{0}$. But then $w_{2c,1}=0$, a contradiction. This also clears reduction type IV. We leave the last three remaining cases to the reader. 
\end{proof}

\subsection{The edge lengths}\label{ProofMainThm3}

We now prove Theorem \ref{MainThm3}.

\subsubsection{Case I}

There are no non-trivial edges in this case.

\subsubsection{Case II}

In terms of the universal family given in \ref{CaseII}, the edge length in the tree is given by $v(t)$ and we have $v(J_{10})=2v(t)$. Since the edge length is doubled, we find that the cycle has length $v(J_{10})$. We moreover have $v(I_{12})=0$, so the desired formula follows. 

\subsubsection{Case III}

We calculate $I_{12}$ and $J_{10}$ for the universal family in \ref{CaseIII}. $J_{10}$ is divisible by $t_{1}^2t_{2}^2$ and the reduction of $m_{10}=J_{10}/t_{1}^2t_{2}^2$ modulo $(t_{1},t_{2})$ is 
$$\overline{m}_{10}=1/4096\cdot{}\lambda_{1}^2\lambda_{2}^4\lambda_{3}^2(\lambda_{2}-1)^4,$$ which is nonzero by our assumptions. The invariant $I_{12}$ is not divisible by either $t_{i}$, but it is zero modulo $(t_{1},t_{2})$. Calculating modulo monomials in the $t_{i}$ of total degree $\geq{3}$, we find that 
\begin{equation}
I_{12}\equiv{1/256\cdot{}(\lambda_{2}^6(\lambda_{2}-1)^4\lambda_{1}^2t_{1}^2+\lambda_{3}^2\lambda_{2}^4(\lambda_{2}-1)^6t_{2}^2)}.
\end{equation}
We can now switch $0$ and $1$ using a projective transformation and assume that $v(t_{1})\leq{v(t_{2})}$. If $v(t_{1})<v(t_{2})$, then $v(I_{12})=2v(t_{1})$ and $2v(t_{2})=v(J_{10})-2v(t_{1})$. If $v(t_{1})=v(t_{2})$, then $v(t_{1})=v(J_{10})/4$. This gives the formulae in the theorem.

\subsubsection{Case IV}\label{CaseIVEdges}

The $v(t_{i})$ from the universal family in \ref{CaseIV} can come in different flavors. By applying another M\"{o}bius transformation, 
we can however assume that $v(t_{1})\leq{v(t_{2})}\leq{v(t_{3})}$. There are then four cases:
\begin{align*}
v(t_{1})<v(t_{2})<v(t_{3}),\\
v(t_{1})=v(t_{2})<v(t_{3}),\\
v(t_{1})<v(t_{2})=v(t_{3}),\\
v(t_{1})=v(t_{2})=v(t_{3}).
\end{align*}
We first calculate $J_{4}$, $I_{12}$ and $J_{10}$. We have $J_{4}=t_{1}^4m_{4}$,  $I_{12}=t_{1}^{12}m_{12}$ and $J_{10}=t_{1}^{12}t_{2}^{2}t_{3}^{2}m_{10}$. Here the reduction of $m_{10}$ modulo $(t_{1},t_{2},t_{3})$ is $1/{4096}\cdot{}\lambda_{1}^2\lambda_{2}^{8}(\lambda_{3,1}-\lambda_{3,2})^2$, which is always nonzero. 
We now calculate $m_{4}$ modulo monomials in the $t_{i}$ of total degree $\geq{3}$ and find 
\begin{equation}
m_{4}\equiv{}1/16\cdot(\lambda_{2}^6t_{1}^2+\lambda_{2}^2(\lambda_{3,1}-\lambda_{3,2})^2t_{2}^2+\lambda_{1}^2\lambda_{2}^2t_{3}^2).
\end{equation}

We then similarly calculate $m_{12}$ modulo the monomial ideal $$J=(t_{1}^3t_{2}^2,t_{1}^3t_{3}^2,t_{2}^2t_{3}^3,t_{2}^3t_{3}^2,t_{1}^2t_{3}^3,t_{1}^2t_{2}t_{3}^2,t_{1}t_{2}^2t_{3}^2,t_{1}^2t_{2}^2t_{3},t_{1}^2t_{2}^3)$$
and find

\begin{equation}
m_{12}\equiv {}1/{256}\cdot{}( \lambda_{2}^{12}(\lambda_{3,1}-\lambda_{3,2})^2t_{1}^{2}t_{2}^2+\lambda_{1}^2\lambda_{2}^{12}t_{1}^2t_{3}^2+\lambda_{1}^2\lambda_{2}^{8}(\lambda_{3,1}-\lambda_{3,2})^2t_{2}^2t_{3}^2).
\end{equation}
Note that the monomials in $J$ are those whose "valuations" are greater than $2v(t_{1})+2v(t_{2})$, $2v(t_{1})+2v(t_{3})$ or $2v(t_{2})+2v(t_{3})$. 

In the first of the four cases introduced above, we find that $v(m_{4})=2v(t_{1})$ and thus $v(J_{4})=6v(t_{1})$. We then also find $v(m_{12})=2v(t_{1})+2v(t_{2})$. This gives us $2v(t_{3})=v(J_{10})-2v(t_{2})-12v(t_{1})$, so we are able to express the $v(t_{i})$ in terms of the tropical Igusa invariants. For the second case, we find that $v(m_{12})=2v(t_{1})+2v(t_{2})=4v(t_{1})$. Using the formula for $J_{10}$, we then also recover $v(t_{3})$. In the third case we again have $v(m_{4})=2v(t_{1})$ and one recovers $v(t_{2})=v(t_{3})$ from $J_{10}$. Finally, if all the $v(t_{i})$ are equal, then we can find them using the formula 
$v(J_{10})=16v(t_{i})$. One now easily verifies that the formulae in the theorem hold.

\subsubsection{Case V}

Recall from \ref{CaseV} that $\overline{J}_{2}\neq{0}$ and $v(J_{10})=6v(t_{1})$. This directly gives the formula in the theorem. 

\subsubsection{Case VI}

We first calculate $\overline{J}_{2}=3/4$. We then recover $v(t_{1})$ from $v(I_{12})$ as in Section \ref{CaseVI} and $v(t_{2})$ from $v(J_{10})$. The formula in the theorem follows directly from this.

\subsubsection{Case VII}

By applying a M\"{o}bius transformation fixing $\infty$, we can assume that $v(t_{2})\leq{v(t_{3})}$. The reduction of $J_{2}$ modulo $(t_{1},t_{2},t_{3})$ is $3/4$, so it is nonzero. By our calculations in Section \ref{CaseVII}, we then have $v(I_{4})=2v(t_{1})$, so this gives the length of the connecting edge. We now calculate the reduction of $m_{12}=I_{12}/t_{1}^6$ modulo the ideal $J=(t_{2}^3,t_{3}^3,t_{1}t_{2}^2,t_{1}t_{3}^2,t_{2}t_{3}^2,t_{3}t_{2}^2)$ and obtain 
\begin{equation}
m_{12}\equiv{} 1/256\lambda_{1}^4(\lambda_{1}^2\lambda_{3}^2t_{3}^2+(\lambda_{2,1}-\lambda_{2,2})^2t_{2}^2).
\end{equation}  
If $v(t_{2})<v(t_{3})$, then we find that $v(m_{12})=2v(t_{2})$ and using $v(J_{10})$ we then also find $v(t_{3})$. If $v(t_{2})=v(t_{3})$, then we can use $v(J_{10})$ to obtain both $v(t_{2})$ and $v(t_{3})$, as before.

\center
\bibliographystyle{alpha}
\bibliography{bibfiles}{}
\end{document}